\documentclass[11pt]{article}
\usepackage{amsmath,amssymb,mathrsfs}
\usepackage{amsfonts}
\input xy
\xyoption{all}

\newtheorem{Thm}{Theorem}[section]
\newtheorem{Prop}[Thm]{Proposition}
\newtheorem{Lemma}[Thm]{Lemma}
\newtheorem{Cor}[Thm]{Corollary}

\newcommand{\pf}{\noindent{\em Proof.}\ }
\newcommand{\qed}{\hfill $\Box$\\}

\renewcommand{\ni}{\noindent}

\newcommand{\ZZ}{\mathbb{Z}}
\newcommand{\QQ}{\mathbb{Q}}
\newcommand{\CC}{\mathbb{C}}
\newcommand{\FF}{\mathbb{F}}
\newcommand{\PP}{\mathbb{P}}
\renewcommand{\AA}{\mathbb{A}}
\newcommand{\GG}{\mathbb{G}}
\newcommand{\tensor}{\otimes}
\newcommand{\uc}{\cup}

\newcommand{\rc}{\subset}

\renewcommand{\projlim}{\underleftarrow{\lim}}

\newcommand{\Spec}{\operatorname{Spec}}

\newcommand{\poch}[2]{\left({#1}\right)_{#2}}
\newcommand{\pochf}[3]{\left(\frac{#1}{#2}\right)_{#3}}
\newcommand{\hyperg}[9]{{}_{#1}F_{#2}\left(\begin{array}{c}
	{#3},{#4},\cdots,{#5} \\
	{#6},{#7},\cdots,{#8}\end{array};{#9}\right)}
\newcommand{\shyperg}[1]{{}_{n}F_{n-1}\left(\begin{array}{c}
	{\frac{1}{n+1}},{\frac{2}{n+1}},\cdots,{\frac{n}{n+1}} \\
	{1},{1},\cdots,{1}\end{array};{#1}\right)}

\newcommand{\lam}{\lambda}
\newcommand{\ve}{\varepsilon}
\newcommand{\om}{\omega}
\newcommand{\HH}{\mathcal{H}}
\newcommand{\MM}{\mathcal{M}}
\newcommand{\RR}{\mathcal{R}}
\newcommand{\TT}{\mathcal{T}}
\newcommand{\VV}{\mathcal{V}}
\newcommand{\Dw}{\mathcal{P}}
\newcommand{\BR}{\mathcal{A}}

\begin{document}

\title{Variation of the Unit Roots
	along the Dwork Family of Calabi-Yau Varieties
	\footnote{This work is supported in part by Professor N. Yui's
		Discovery Grant from NSERC, Canada.}}
\author{\sc Jeng-Daw Yu\footnote{
	email address: {\tt jdyu@mast.queensu.ca}}}
\date{August 2007}
\maketitle

\begin{abstract}
We study the variation of unit roots
of the Dwork families of Calabi-Yau varieties over a finite field
by the method of Dwork-Katz
and also from the point of view of formal group laws.
A $p$-adic analytic formula
for the unit roots
away from the Hasse locus
is obtained.
\end{abstract}

\section{Introduction}

Let us first look at an example
and then explain the main results of this paper.\\

\ni{\em (a) The Legendre family}\\

Recall the following classical results
(see \cite{Dwork}, \S $6(i)$ or \cite{Katz_Dwork}, \S 8).\\

Consider the Legendre family of elliptic curves $E_{\lam}$
whose affine part is given by
\[ E_{\lam}: y^2 = x(x-1)(x - \lam) \]
with the parameter $\lam \neq 0, 1$.
Over the complex numbers $\lam \in \CC$,
the relative de Rham cohomology $H^1_{dR}$ of degree 1
of the family is free of rank 2.
The Hodge filtration ${\rm Fil}^1 \rc H^1_{dR}$
is generated by the differential of the first kind
\[ \om = \frac{dx}{y}. \]
Let $\nabla$ be the Gauss-Manin connection on $H^1_{dR}$.
Then $\om$ satisfies the associated Picard-Fuchs equation
$\nabla(L) \om = 0$,
where
\begin{equation}\label{PF_Legendre}
L = \lam (\lam - 1) \frac{d^2}{d\lam^2}
	+ (1 - 2\lam) \frac{d}{d\lam} - \frac{1}{4}.
\end{equation}
Up to a constant,
the unique holomorphic solution to
the differential equation (\ref{PF_Legendre})
at $\lam = 0$
is given by the Gauss hypergeometric series
\begin{equation}\label{F21}
F(\lam) = {}_2F_1 \left( \frac{1}{2}, \frac{1}{2}; 1; \lam \right).
\end{equation}
Furthermore,
the section
\begin{equation}\label{hs_Legendre}
u = \lam(1-\lam) F \om' - \lam(1-\lam) F' \om
\end{equation}
of $H^1_{dR}$
is a global horizontal section
with respect to $\nabla(d/d\lam)$.
Here $\om' = \nabla(d/d\lam)\om$
and $F' = dF/d\lam$
are the derivatives.\\

Interestingly
the series (\ref{F21}) also gives information
for the Legendre family over a finite field.
More precisely,
let $p$ be an odd prime.
Let $H(x) = F(x)^{<p}$
be the truncation of $F(x)$
up to degree $p-1$.
Let $\lam \in \FF_q, \lam \neq 0, 1$,
where $q = p^s$.
Then the elliptic curve $E_{\lam}$
is ordinary
if and only if $H(\lam) \neq 0$ in $\FF_q$.\\

Moreover,
let $W(\FF_q)$ be the ring of Witt vectors of $\FF_q$.
Let $\hat{\lam} \in W(\FF_q)$
be the Teichm\"uller lifting of $\lam$.
Suppose $H(\lam) \neq 0$.
Then the formal power series
\[ f(x) = \frac{F(x)}{F(x^p)} \]
converges at $\hat{\lam}$ as a series over $W(\FF_q)$.
Write the zeta function of $E_{\lam}$ over $\FF_q$ as
\[ Z(E_{\lam},T) = \frac{1 - aT + qT^2}{(1-T)(1-qT)}. \]
Then $a = \pi + \pi'$ with $\pi \pi' = q$ and
\begin{equation}\label{Leg_ur}
\pi = \ve^s f(\hat{\lam}) f(\hat{\lam}^2) \cdots f(\hat{\lam}^{p^{s-1}})
\end{equation}
with $\ve = (-1)^{(p-1)/2}$.
The algebraic integer $\pi$ is a $p$-adic unit
and is called the {\em unit root}
of the ordinary elliptic curve $E_{\lam}$.\\

We remark that
the formula (\ref{Leg_ur}) can be derived easily
by the method in \S \ref{SB},
based on formal group laws.\\

\ni{\em (b) The Dwork family}\\

In this paper,
we shall generalize the above results
to certain higher dimensional cases.\\

Throughout the paper,
we let $n \geq 2$ be an integer.
Let $\BR$ be a base ring
in which $(n+1)$ is invertible.
Set
\[ X = [X_1: X_2: \dots : X_{n+1}] \]
to be the homogeneous coordinates
of the projective space $\PP^n$ over $\BR$.
We will write $t \in \AA^1 \uc \{ \infty \}$
as the coordinate of $\PP^1$.\\

\ni{\em Definition}
(cf. \cite{Katz_AL}, \S 1).
The {\em Dwork family} over $\BR$
is the one-parameter family $V_t$
of Calabi-Yau hypersurfaces in $\PP^n$ over $t \in \PP^1$
defined by the equation $\Dw_t(X) = 0$,
where
\begin{equation}\label{Dwork}
\Dw_t(X) = X_1^{n+1} + X_2^{n+1} + \cdots + X_{n+1}^{n+1}
	- (n+1) t X_1 X_2 \cdots X_{n+1}.
\end{equation}
We also set $\VV$ to be the total space
of the family in $\PP^n \times \PP^1$
and $v: \VV \to \PP^1$ to be the fiber map.\\

In this paper,
we study the variation along $t$
of the zeta function of the Dwork family $V_t$
over a finite field $\BR = \FF_q$.
Relevant definitions will be given in \S \ref{Prel}.
As in the case of Legendre family,
we show that the zeta function of $V_t$
is closely related to the unique holomorphic solution $F$
of the associated Picard-Fuchs equation of the family.
For the Dwork family,
it turns out that $F$ is a generalized hypergeometric series.
We shall define the Hasse invariant for the Dwork family.
Similar to the formula (\ref{Leg_ur}),
we shall derive a formula for the unit root of $V_t$
in terms of
the ratio, $f$, of $F$ and its Frobenius twist
(Theorem \ref{unitroot}).
We drive the formula by two different methods,
(I) and (II),
which we describe briefly now.
\renewcommand{\labelenumi}{{\rm (\Roman{enumi})}}
\begin{enumerate}
\item
We follow Katz's crystalline interpretation \cite{Katz_Dwork}
of Dwork's work on the variation of zeta functions
of hypersurfaces \cite{Dwork}.
We generalize the construction
of horizontal sections (\ref{hs_Legendre})
to the case of the Dwork family (Corollary \ref{HS_Dwork}).
In this way,
we see immediately
that the Picard-Fuchs equation enters the picture.

\item
We study the Frobenius action via formal groups.
We determine the formal group
associated to $H^{n-1}(V_t, W\mathcal{O})$ of $V_t$ explicitly
by writing down a formal group law $G_t$ for it
(Proposition \ref{law}).
For this, we will follow the work of Stienstra \cite{S}.
From this point of view,
the series $F$ appears
as a certain $p$-adic limit of the coefficients
of the logarithm of $G_t$.
\end{enumerate}

Besides the formula for the unit root,
each method also provides bonus information
of different flavor.
Method (I) shows that
$F_i/F$ has a $p$-adic analytic continuation,
where $F_i$ is the $i$-th derivative of $F$.
Method (II) provides a good approximation
to the unit root similar to \cite{Dwork}, Lemma (3.4)(i).
The second method might be viewed
as a dual approach of \cite{Katz_Internal}.
In op.cit.,
the author studied the highest slope part
while here we look at the slope zero part directly.\\

We study the relation
between the Picard-Fuchs equation
and the unit root
in a hope that,
at least over positive characteristic,
one can study the arithmetic of the Dwork family
or other families of Calabi-Yau varieties
inductively
by reducing the weights
of the cohomology.
For example,
in the threefold case ($n = 4$),
the Picard-Fuchs equation studied here
takes care of the most transcendental part
of the middle cohomology $H^3$.
Then the remaining part can be viewed
as a family of abelian varieties
by the construction of intermediate Jacobians.
One might combine these two pieces
to obtain a better understanding
of the whole $H^3$.
Also
the methods developed in this paper
should be extendable
to families of Calabi-Yau varieties
of generalized hypergeometric type
(e.g.
complete intersections
in weighted projective spaces).\\

The original motivation of studying the associated formal group
for the Dwork family
was to see if one can get a geometric interpretation
of the congruence \cite{Dwork}, Lemma (3.4)(i),
and eventually the congruences in \S 1, Corollary 2 therein.
If this is the case,
one might be able to prove the similar congruences
for different families of Calabi-Yau varieties
of non-hypergeometric type.
Numerical computation suggests that
the {\em Ap\'ery numbers} listed in \cite{SB}, Table 7
satisfy the same type of congruences in loc.cit.
Unfortunately,
Method (II) gives rise to a family
of hypergeometric functions
different from the truncations of $F$,
and the congruences for Ap\'ery numbers
remain as an open question.\\

Recently,
there are several papers
dealing with the Dwork family.
In \cite{HST},
the authors compute the Zariski closure
of the monodromy group
in characteristic 0 of this family
and apply it
to the study of the Sato-Tate conjecture.
The papers \cite{Katz_AL} and \cite{RW}
also investigate the family
in characteristic $p$
via $\ell$-adic Fourier transforms.
The moment zeta functions of the Dwork family
are computed in \cite{RW}.
In \cite{Kl},
the zeta functions
for more general monomial deformations
of Fermat type hypersurfaces
in weighted projective spaces
are studied
via Dwork's deformation theory.\\

The author wants to
thank Prof. N. Yui
for her interest in this work,
helpful discussion,
and her encouragement.
She also provided useful suggestions and commends
to the first draft of this paper.
Thanks are also due to Prof. N. Katz
for bringing my attention
to \cite{Katz_Internal}.

\section{Preliminary}\label{Prel}

In this section,
we give relevant definitions
that are used in this paper.\\

Fix an integer $n \geq 2$.
Recall that $\BR$ is the base ring
in which $(n+1)$ is invertible.
The map $v:  \VV \to \PP^1$
indicates the fiber map of the Dwork family
(see (\ref{Dwork}))
\[ V_t: \Dw_t(X) = 0. \]
Let
\[ \mu_{n+1} = \Spec \BR[x]/(x^{n+1} - 1) \]
be the group scheme
of the $(n+1)$-st roots of unity over $\BR$.
Note that
over the subscheme $\TT = \AA^1 \setminus \mu_{n+1}$ in $\PP^1$,
the map $v$ is smooth.
Let $\HH = \mu_{n+1}^{n+1}/\mu_{n+1}$,
where the quotient is via the diagonal embedding.
Let $\HH_0$ be the subgroup of $\HH$ defined by
\[ \HH_0 = \{ \underline{\zeta} = (\zeta_1, \dots, \zeta_{n+1}) \in \HH \ |\
	\zeta_1 \cdots \zeta_{n+1} = 1 \}. \]
Then $\HH_0$ acts on each fiber $V_t$
of the family
by
\[ \underline{\zeta} (x_1, x_2, \dots, x_{n+1})
	= (\zeta_1 x_1, \zeta_2 x_2, \dots, \zeta_{n+1} x_{n+1}) \]
for $x = (x_1, \dots, x_{n+1}) \in V_t$.
We are mainly interested in
the relative cohomology of $v$ of degree $(n-1)$
fixed by $\HH_0$.
We discuss the details in the following two cases.\\

\ni{\em (a) Over $\CC$}\\

Suppose $\BR = \CC$.
Let $\MM_{dR}$ be
the fixed part by $\HH_0$
of the relative de Rham cohomology
$R^{n-1}v_* \Omega^{\bullet}_{\mathcal{V/T}}$ of degree $(n-1)$
over $\TT$.
Here $\Omega^{\bullet}_{\mathcal{V/T}}$
is the complex of relative differential forms.
Then
the sheaf $\MM_{dR}$ is locally free of rank $n$.
Denote ${\rm Fil}^{\bullet}$ the Hodge filtration of $\MM_{dR}$.
At each point,
$\MM_{dR}$ has Hodge numbers
\begin{equation}\label{hij}
h^{i,n-1-i} := \dim_{\CC} {\rm Fil}^i/{\rm Fil}^{i+1} = 1
\end{equation}
for all $0 \leq i \leq n-1$ (\cite{HST}, \S 1).
Let $\nabla$ be the Gauss-Manin connection.
Then $\MM_{dR}$ is stable
under $\nabla$.\\

Let $\Omega$ be the differential $n$-form on $\PP^n$
defined by
\[ \Omega = \sum_{i=1}^{n+1} (-1)^i X_i dX_1 \wedge
	\cdots \wedge \widehat{dX_i} \wedge \cdots dX_{n+1}. \]
Here $\widehat{dX_i}$ means
the deletion of the $i$-th component $dX_i$.
Let
\[ \xi = \text{Res}_{V_t} (\Omega/\Dw_t) \]
be the cohomology class in the top forms
\[ H^0(V_t/\CC, \Omega^{n-1}) \rc H^{n-1}_{dR}(V_t/\CC) \]
given by the residue of the meromorphic differential $n$-form
$\Omega/\Dw_t$.
Let
\[ \xi_i = \left( \nabla \left( \frac{d}{dt} \right) \right)^i \xi \]
be the $i$-th derivative of $\xi$.
It is shown (loc.cit) that
\[ \xi_i \in {\rm Fil}^{n-1-i} - {\rm Fil}^{n-i} \]
and the set $\{ \xi_i \}_{0 \leq i < n}$
generates $\MM_{dR}$.
Finally define a new section $\eta$
of $\MM_{dR}$ by
\begin{equation}\label{eta}
\eta = t \cdot \xi.
\end{equation}
\\

\begin{Prop}
In $\MM_{dR}$,
we have
\renewcommand{\labelenumi}{{\rm (\arabic{enumi})}}
\begin{enumerate}
\item
$\xi$ satisfies the Picard-Fuchs equation
$\nabla(L_t) \xi = 0$,
where
\begin{eqnarray}\label{PF_t}
L_t = \prod_{i = 1}^n \left(\delta + 1 - i \right)
		- t^{n+1} \left(\delta + 1 \right)^n
\end{eqnarray}
with $\delta = \delta_t = t \frac{d}{dt}$.
\item
$\eta$ satisfies the differential equation
$\nabla(L_{\lam}) \eta = 0$, where
\begin{eqnarray}\label{PF}
L_{\lam} = \theta^n
	- \lam \prod_{i = 1}^n \left(\theta + \frac{i}{n+1}\right)
\end{eqnarray}
with $\lam = t^{-(n+1)}$ and
$\theta = \theta_{\lam} = \lam \frac{d}{d\lam}$.
\end{enumerate}
\end{Prop}

\pf
We prove
that (\ref{PF_t}) and (\ref{PF})
are indeed the Picard-Fuchs equations
for $\xi$ and $\eta$, respectively,
as a follow-up of \cite{HST}, \S 1.
To shorten the notation,
we simply write $D \om$
to indicate $\nabla(D) \om$
for any derivative $D$
and de Rham cohomology class $\om$.

We employ the same notations as in \cite{HST}.
Let $\om_1 = \eta = t \cdot \xi$
and define $\om_r$ inductively by
\[ \om_{r+1} = \delta \om_r - \om_r. \]
For all $0 \leq i < j$,
define the constants $A_{i,j}$
by
\[ T^r = \sum_{i=0}^r A_{i,r+1} (T-1)(T-2) \cdots (T-r+i) \]
as polynomials in $T$ for all $r \geq 0$.
Then in the first half of \S 1 in op.cit.
(formula in page 9),
we have, as cohomology classes,
\begin{equation}\label{HST}
\om_{n+1} = t^{n+1} (A_{0,n+1} \om_{n+1} + A_{1,n+1} \om_n
	+ \cdots + A_{n,n+1} \om_1).
\end{equation}

Inductively, we can derive
\begin{equation}\label{rto1}
w_{r+1} = (\delta-1)(\delta-2) \cdots (\delta-r) \om_1,
\end{equation}
and
\begin{equation}\label{1toxi}
\delta^r \om_1 = t(\delta + 1)^r \xi.
\end{equation}
Plugging (\ref{rto1}) into (\ref{HST}),
we get
\begin{eqnarray}
\nonumber
\left( \prod_{i=1}^n (\delta-i) \right) \om_1
	&=& \left( t^{n+1} \sum_{i=0}^n A_{i,n+1}
		(\delta-1)(\delta-2) \cdots (\delta-n+i) \right) \om_1 \\
	&=& t^{n+1} \delta^n \om_1.
	\label{PF_HST}
\end{eqnarray}
Equations (\ref{1toxi}) and (\ref{PF_HST}) imply
\[ t \prod_{i=1}^n (\delta + 1 - i)\xi = t^{n+2} (\delta + 1)^n \xi, \]
which gives the claimed equation (\ref{PF_t}).\\

On the other hand,
set $\lam = t^{-(n+1)}$
and $\theta = \lam \frac{d}{d \lam}$.
Then
\[ \delta = -(n+1) \theta. \]
Put it into (\ref{PF_HST}),
we obtain
\[ \prod_{i=1}^n \left( -(n+1) \theta -i \right) \eta
	= \lam^{-1} \left( -(n+1) \theta \right)^n \eta, \]
which is equivalent to (\ref{PF}).
\qed

We remark that
the discussion here
is valid also to any field
that can be embedded into $\CC$.\\

\ni{\em (b) Over positive characteristic}\\

For relevant descriptions
of crystals over a smooth base,
see \cite{Katz_Dwork}, \S\S 1, 5, 7
and \cite{Katz_Slope}, \S\S 2.1, 2.4.\\

Let $\BR = k$ be a field of characteristic $p > 0$
with $(p, n+1) = 1$.
Let $W = W(k)$ be the ring of Witt vectors of $k$.
For $t \in k, t^{n+1} \neq 1$,
the crystalline cohomology $H^{n-1}_{cris} (V_t/W)$
is a free $W$-module
equipped with an absolute Frobenius action $\phi$.
We will simply call the Newton polygon of $V_t$
to mean the Newton polygon of $H^{n-1}_{cris} (V_t/W)$,
and similarly for the Hodge polygon
(see \cite{Katz_Slope}, \S\S 1.2 and 1.3).\\

\ni{\em Definition.}
We say that $V_t$ is {\em ordinary}
if the Newton polygon coincides the Hodge polygon
of $V_t$
(cf. \cite{Illusie}, \S\S 1.1, 1.3).\\

\begin{Thm}
Let $k$ be a field of characteristic $p > 0$
with $(p, n+1) = 1$.
Then the Dwork family defined by equation (\ref{Dwork})
over $k$ is generically ordinary.
\end{Thm}

\pf
Let $k[[t^{-1}]]$ be the localization
of the parameter space $\PP^1$ at $t = \infty$.
Let $\TT' = \Spec k[[t^{-1}]]$
and $\VV'$ be the restriction of the family $\VV$ to $\TT'$.
The special fiber $V_{\infty}$
is the union of coordinate hyperplanes.
Each arbitrary intersection among them
is isomorphic to some projective space
and is obviously ordinary.
Thus
the generic fiber $\widetilde{\VV'}$ of $\VV'$
is ordinary
(\cite{Illusie}, Proposition 1.10).
\qed

Let $\hat{t} \in W$ be a lifting of $t \in k$.
Then $H^{n-1}_{dR}(V_{\hat{t}}/W)$
is canonically isomorphism to $H^{n-1}_{cris}(V_t/W)$
and the identification is compatible
with the $\HH_0$-action.
Let $\MM_{cris}$ be the fixed part
of $H^{n-1}_{cris}(V_t/W)$  by $\HH_0$.
Then $\MM_{cris}$ is a direct summand
of $H^{n-1}_{cris}(V_t/W)$,
and similarly for the fixed part $\MM_{dR}$
of $H^{n-1}_{dR}(V_{\hat{t}}/W)$ by $\HH_0$.
Since the Newton polygon is on or above the Hodge polygon
of $V_t$,
and by a glance at
the Hodge polygon of $\MM_{dR}$
described in (\ref{hij}),
we see that
the first slope of the Newton polygon of $V_t$
must be $0, 1/2$, or $\geq 1$.\\

Assume now that $k$ is perfect.
Let $\mathfrak{Art}_k$ be the category
of Artinian local $k$-algebras,
and $\mathfrak{AG}$ the category
of abstract abelian groups.
The Artin-Mazur functor
\[ G_t: \mathfrak{Art}_k \to \mathfrak{AG} \]
is defined by
\[ G_t(R, \mathfrak{m}) = {\rm Ker}
	\left\{ H^{n-1}_{et}\left(X \tensor_k R, \hat{\GG}_m\right)  \to
	H^{n-1}_{et}\left(X \tensor_k R/\mathfrak{m}, \hat{\GG}_m\right)
		\right\}, \]
for $(R,\mathfrak{m})$ an object in $\mathfrak{Art}_k$.
The functor $G_t$ is pro-representable
by a one-dimensional commutative formal group
(\cite{S}, Theorem 1).
We call $G_t$ the
{\em formal group associated to $V_t$}.
The (covariant) Cartier module
(of $p$-typical curves) of $G_t$
is canonically isomorphic to $H^{n-1}(V_t, W\mathcal{O})$
as a $W[\phi]$-module
(\cite{Illusie_Witt}, Remarque II.2.15).
By the description of the Newton polygon
in the last paragraph,
$G_t$ is of height $1, 2$ or $\infty$.
Notice that
the formal group $G_t$
can be defined over a more general base ring
$\BR$ (see \cite{S}, Theorem 1).\\

Finally suppose $k = \FF_q$ is a finite field
of $q$ elements.
Suppose $t \in \FF_q, t^{n+1} \neq 1$.
If the first slope of the Newton polygon of $V_t$ is zero,
there is a unique $p$-adic unit root
of the geometric Frobenius endomorphism
acting on $H^{n-1}_{cris}(V_t/W)$.
We will call this the {\em unit root} of $V_t$.\\

\section{Existence of a global horizontal section}\label{Sec_HS}

Here we explicitly construct a global horizontal section
for the Dwork family
with respect to the Gauss-Manin connection
over characteristic zero.\\

\begin{Lemma}\label{ai=0}
For any two positive integers $k, m$ with $k-1 < m$,
we have
\[ \sum_{r = 0}^m (-1)^r \binom{m-k+r}{k-1} \binom{m}{r} = 0. \]
\end{Lemma}

\pf
Consider the function
$a(x) = x^{m-k} (1+x)^m$.
Then
\[ (k-1)! \sum_{r = 0}^m (-1)^r \binom{m-k+r}{k-1} \binom{m}{r}
	= \pm \frac{d^{k-1}}{dx^{k-1}} a(-1) = 0 \]
for $k-1 < m$.
\qed

\begin{Lemma}\label{neven}
Consider functions
$b^{(m)}_i = b^{(m)}_i (a_1, \dots, a_m), 0 \leq i \leq 2m$,
of $m$ variables $a_1, \dots, a_m$ defined by
\[ b(x) = \prod_{i = 1}^m (x + a_i) (x + 1 - a_i)
	= \sum_{i = 0}^{2m} b^{(m)}_i x^{2m - i}. \]
Then for $1 \leq k \leq m$,
we have
\[ \sum_{i=0}^{2m-2k+1} (-1)^i \binom{k-1+i}{k-1} b^{(m)}_{2m-2k+1-i} = 0. \]
\end{Lemma}

\pf
We prove this by deformation of $a_i$ and by induction on $m$.
It is easy to establish the equality when $m = 1$.
Then notice first that
\begin{eqnarray*}
\frac{\partial b(x)}{\partial a_i} &=& \frac{b(x)}{x+a_i} - \frac{b(x)}{x+1-a_i} \\
	&=& (1-2a_i) \frac{b(x)}{(x+a_i)(x+1-a_i)}.
\end{eqnarray*}
Comparing the coefficients on both sides, we have
\[ \frac{\partial}{\partial a_i} b^{(m)}_k (a_1, \dots, a_m)
	= (1-2a_i) b^{(m-1)}_{k-2} (a_1, \dots, \hat{a_i}, \dots, a_m). \]
Here $\hat{a_i}$ means the deletion
of the $i$-th component $a_i$.
Thus by induction,
\begin{eqnarray*}
&& \frac{\partial}{\partial a_i} \sum_{i=0}^{2m-2k+1} (-1)^i
		\binom{k-1+i}{k-1} b^{(m)}_{2m-2k+1-i} \\
&=& (1-2a_i) \sum_{i=0}^{2m-2k-1} (-1)^i
		\binom{k-1+i}{k-1} b^{(m-1)}_{2m-2k-1-i} \\
&=& 0.
\end{eqnarray*}
Therefore $\sum (-1)^i \binom{k-1+i}{k-1} b^{(m)}_{2m-2k+1-i}$
is a constant.

Secondly for $a_1 = \dots = a_m = 0$,
we have $b^{(m)}_i = \binom{m}{i}$ and then
\begin{eqnarray*}
\sum_{i=0}^{2m-2k+1} (-1)^i \binom{k-1+i}{k-1} b^{(m)}_{2m-2k+1-i}
	&=& \sum_{i = 0}^{2m-2k+1} (-1)^i
			\binom{k-1+i}{k-1} \binom{m}{2m-2k+1-i} \\
	&=& \sum_{i = m-2k+1}^{2m-2k+1} (-1)^i
			\binom{k-1+i}{k-1} \binom{m}{2m-2k+1-i} \\
	&=& \sum_{r = 0}^m (-1)^r \binom{m-k+r}{k-1} \binom{m}{r} \\
	&=& 0
\end{eqnarray*}
by Lemma \ref{ai=0} and this completes the proof.
\qed

\begin{Lemma}\label{nodd}
Consider functions
$\beta_i = \beta_i (a_1, \dots, a_m), 0 \leq i \leq 2m+1$,
of variables $a_1, \dots, a_m$ defined by
\[ \beta(x) = \left(x+\frac{1}{2}\right) \prod_{i = 1}^m (x + a_i) (x + 1 - a_i)
	= \sum_{i = 0}^{2m+1} \beta_i x^{2m +1 - i}. \]
Then
\[ \sum_{i=0}^{2m} (-1)^i \beta_{2m-i} = 2 \beta_{2m+1}, \]
and for $1 \leq k \leq m$,
we have
\[ \sum_{i=0}^{2m-2k+1} (-1)^i \binom{k-1+i}{k-1} \beta_{2m-2k+1-i}
	=\frac{1}{2} \sum_{i=0}^{2m-2k} (-1)^i \binom{k+i}{k} \beta_{2m-2k-i}. \]
\end{Lemma}

\pf
We define $b_i = b^{(m)}_i (a_1, \dots, a_m), 0 \leq i \leq 2m$
as in Lemma \ref{neven}.
Let $b_{-1} = b_{2m+1} = 0$.
Then
\[ \beta_i = b_i + \frac{1}{2} b_{i-1} \]
for $0 \leq i \leq 2m+1$.
Thus by Lemma \ref{neven},
\begin{eqnarray*}
\sum_{i=0}^{2m} (-1)^i \beta_{2m-i}
	&=& \sum_{i=0}^{2m}
		(-1)^i \left( b_{2m-i} + \frac{1}{2} b_{2m-i-1} \right) \\
	&=& b_{2m} - \frac{1}{2} \sum_{i=0}^{2m-1} (-1)^i b_{2m-1-i} \\
	&=& b_{2m} \\
	&=& 2 \beta_{2m+1}.
\end{eqnarray*}
On the other hand,
\begin{eqnarray*}
&& \sum_{i=0}^{2m-2k+1} (-1)^i \binom{k-1+i}{k-1} \beta_{2m-2k+1-i} \\
	&=& \sum_{i=0}^{2m-2k+1} (-1)^i \binom{k-1+i}{k-1}
		\left(b_{2m-2k+1-i} + \frac{1}{2} b_{2m-2k-i} \right) \\
	&=& \sum_{i=0}^{2m-2k+1} (-1)^i \binom{k-1+i}{k-1} b_{2m-2k+1-i}
		+ \frac{1}{2} \sum_{i=0}^{2m-2k+1}
			(-1)^i \binom{k-1+i}{k-1} b_{2m-2k-i} \\
	&=&\frac{1}{2} \sum_{i=0}^{2m-2k} (-1)^i \binom{k-1+i}{k-1} b_{2m-2k-i}
\end{eqnarray*}
by Lemma \ref{neven}, and since $b_{-1} = 0$.
Similarly,
\begin{eqnarray*}
\sum_{i=0}^{2m-2k} (-1)^i \binom{k+i}{k} \beta_{2m-2k-i}
	= \sum_{i=0}^{2m-2k} (-1)^i \binom{k+i}{k}
		\left(b_{2m-2k-i} + \frac{1}{2} b_{2m-2k-1-i} \right),
\end{eqnarray*}
and
\[ \sum_{i=0}^{2m-2k} (-1)^i \binom{k+i}{k} b_{2m-2k-1-i} = 0. \]
Thus we have the desired equality.
\qed

\begin{Thm}\label{HS}
Let $A$ be a ring.
Let $B$ and $M$ be two $A[x]/A$-differential modules,
where $x$ is a variable.
Fix an $A[x]/A$-differential $D$
and assume $e \in A[x]$ satisfying $De = e$.
Suppose $\{a_i\} \rc A$ is stable
under the transformation $a_i \mapsto 1 -a_i$.
Suppose that
$g \in B$ and $\eta \in M$ satisfy the differential equation
$Lv = 0$, where
\begin{equation}
L = D^{n} - e \prod_{i = 1}^n (D + a_i)
	= D^n - e \sum_{i=0}^n b_i D^{n-i}.
\end{equation} 
Write $n = 2m - \ve$, where $m \in \ZZ, \ve = 0$ or $1$.
Set $g^{(i)} = D^i g$, and $\eta^{(j)} = D^j \eta$.
Let
\[ c_{ij} = \sum_{r=0}^j (-1)^r \binom{i+r-1}{i-1} b_{j-r}. \]
Then the element $u \in B \tensor_{A[x]} M$ defined by
\begin{eqnarray}
\nonumber
u = && (1-e)\sum_{i=0}^{n-1} (-1)^i g^{(i)}\eta^{(n-1-i)} \\
\nonumber
	&+& e \sum_{i=1}^{m-1}\sum_{j=1}^{n-2i} (-1)^i c_{ij}
		\left[g^{(i-1)}\eta^{(n-i-j)} - (-1)^{\ve} g^{(n-i-j)}\eta^{(i-1)}\right] \\
	&+& \ve e \sum_{i=1}^{m-1} (-1)^i c_{i, n+1-2i} g^{(i-1)} \eta^{(i-1)}
\end{eqnarray}
is a horizontal section with respect to $D$.
\end{Thm}

\pf
Notice that by an easy calculation,
we have
\[ c_{i+1, 0} = b_0 = 1, \]
\[ c_{1,j} + c_{1,j+1} = b_{j+1}, \]
and
\begin{equation}\label{eqn_hs}
c_{i+2,j} + c_{i+2,j+1} = c_{i+1,j+1}
\end{equation}
for all $i, j \geq 0$.
To simplify the notation,
we let
\[ (g_i, \eta_j) := g^{(i)}\eta^{(j)} - (-1)^{\ve} g^{(j)}\eta^{(i)}. \]
We have
\begin{eqnarray*}
Du = && (1-e) (g, \eta_n) - e\sum_{i=0}^{m-1} (-1)^i (g_i, \eta_{n-1-i}) \\
	&+& e \sum_{i=1}^{m-1}\sum_{j=1}^{n-2i} (-1)^i c_{ij}
		\Big( (g_{i-1}, \eta_{n+1-i-j}) + (g_i, \eta_{n-i-j})
			+ (g_{i-1}, \eta_{n-i-j}) \Big) \\
	&+& \ve e \sum_{i=1}^{m-1} (-1)^i c_{i, n+1-2i}
		\left( (g_{i-1}, \eta_i) + \frac{1}{2} (g_{i-1}, \eta_{i-1}) \right).
\end{eqnarray*}

We now distinguish two cases.

(a)
Suppose $n$ is even, i.e. $\ve = 0$.
Letting $k = i+1$ in Lemma \ref{neven},
we get
\[ c_{i+1, n-1-2i} = 0 \]
for $0 \leq i \leq m-1$.
To simplify the notation,
we let
\[ [g_i, \eta_j] := g^{(i)}\eta^{(j)} - g^{(j)}\eta^{(i)}. \]
Collecting all terms of the form $[g, \eta_j]$ in $Du$,
we have
\begin{eqnarray*}
&&	(1-e)[g, \eta_n] - e\left\{ (1 + c_{1,1}) [g, \eta_{n-1}]
	+ \sum_{j=2}^{n-2} (c_{1,j-1} + c_{1,j}) [g, \eta_{n-j}]
	+ c_{1,n-2} [g, \eta_1] \right\} \\
&=& (1-e)[g, \eta_n] - e \sum_{j=1}^{n-1} (c_{1,j-1} + c_{1,j}) [g, \eta_{n-j}] \\
&=& (1-e)[g, \eta_n] - e \sum_{j=1}^{n-1} b_j [g, \eta_{n-j}] \\
&=& g [(L + eb_n)\eta] - \eta [(L + eb_n)g] \\
&=& 0.
\end{eqnarray*}
The last equality follows
since $Lg = L\eta = 0$.
On the other hand,
the coefficient of $[g_i, \eta_j]$ for $0 < i < j$ in $Du$ is
given by
\begin{eqnarray*}
&& (-1)^{i+1} e \times \left\{ \begin{array}{lll}
	(1 + c_{i+1,1} - c_{i,1}) &{\rm if}& j = n-1-i \\
	(c_{i+1,n-1-i-j} + c_{i+1,n-i-j} - c_{i,n-i-j}) &{\rm if}& n-2-i \geq j \geq 2+i \\
	(c_{i+1,n-2-2i} - c_{i,n-1-2i}) &{\rm if}& j = 1+i
	\end{array} \right. \\
&=& (-1)^{i+1} e (c_{i+1,j-1} + c_{i+1,j} - c_{i,j}) \\
&=& 0.
\end{eqnarray*}
The last equality follows from (\ref{eqn_hs}).
Thus $Du = 0$ in this case.

(b)
Suppose $n$ is odd, i.e. $\ve = 1$.
Letting $k = i$ in Lemma \ref{nodd},
we get
\[ c_{1,n-1} = 2 b_n \]
and
\[ c_{i,n-2i} = \frac{1}{2} c_{i+1,n-1-2i} \]
for $1 \leq i \leq m-1$.
To simplify the notation,
we let
\[ \{g_i, \eta_j\} := g^{(i)}\eta^{(j)} + g^{(j)}\eta^{(i)}. \]
Collecting all terms of the form $\{g, \eta_j\}$ in $Du$,
we have
\begin{eqnarray*}
&&	(1-e)\{g, \eta_n\} - e\left[ (1 + c_{1,1}) \{g, \eta_{n-1}\}
	+ \sum_{j=2}^{n-1} (c_{1,j-1} + c_{1,j}) \{g, \eta_{n-j}\}
	+ c_{1,n-1} g \eta \right] \\
&=& (1-e)\{g, \eta_n\} - e \left[
	\sum_{j=1}^{n-1} (c_{1,j-1} + c_{1,j}) \{g, \eta_{n-j}\}
	+ b_n \{g, \eta\} \right] \\
&=& (1-e)\{g, \eta_n\} - e \sum_{j=1}^{n} b_j \{g, \eta_{n-j}\} \\
&=& g L\eta + \eta Lg \\
&=& 0.
\end{eqnarray*}
On the other hand,
the coefficient of $\{g_i, \eta_j\}$ for $0 < i \leq j$ in $Du$ is
given by
\begin{eqnarray*}
&& (-1)^{i+1} e \times \left\{ \begin{array}{lll}
	(1 + c_{i+1,1} - c_{i,1}) &{\rm if}& j = n-1-i \\
	(c_{i+1,n-1-i-j} + c_{i+1,n-i-j} - c_{i,n-i-j}) &{\rm if}& n-2-i \geq j \geq 1+i \\
	(\frac{1}{2}c_{i+1,n-1-2i} - c_{i,n-2i}) &{\rm if}& j = i
	\end{array} \right. \\
&=& (-1)^{i+1} e \times \left\{ \begin{array}{lll}
	(c_{i+1,n-1-i-j} + c_{i+1,n-i-j} - c_{i,n-i-j}) &{\rm if}& n-1-i \geq j \geq 1+i \\
	(\frac{1}{2}c_{i+1,n-1-2i} - c_{i,n-2i}) &{\rm if}& j = i
	\end{array} \right. \\
&=& 0.
\end{eqnarray*}
Thus $Du = 0$.
\qed

\ni{\em Remark.}
If we write the element $u$ in Theorem \ref{HS} as
\[ u = \sum_{i=0}^{n-1} C_i \eta^{(n-1-i)}, \]
then $C_0 = (1-e)g$
and $C_i$ is an $A[e]$-linear combination
of $\{ g, g^{(1)}, \dots, g^{(i)} \}$.\\

Now we go back to the Dwork family over $\CC$.
The differential equation (\ref{PF})
has a unique power series solution
$F(\lam) \in \CC[[\lam]]$ with constant term 1,
which is holomorphic near $\lam = 0$.
We know explicitly that
$F(\lam)$ is given by
a hypergeometric series:
\begin{eqnarray}\label{F_n_n-1}
F(\lam) &=& \shyperg{\lam} \\
	\nonumber
	&=& \sum_{r=0}^{\infty} \frac{\pochf{1}{n+1}{r}
		\pochf{2}{n+1}{r} \cdots\pochf{n}{n+1}{r}}{(r!)^n}
		\lam^r,
\end{eqnarray}
where $\poch{a}{0} = 1$
and $\poch{a}{r} = a(a+1) \cdots (a+r-1)$ for $r > 0$
is the {\em Pochhammer symbol}.\\

\begin{Cor}\label{HS_Dwork}
Consider the Dwork family
over $\BR = \CC$.
Let $g = F(\lam)$ be the hypergeometric series (\ref{F_n_n-1}).
Let $\eta$ be the de Rham cohomology class
given by (\ref{eta}).
Then up to a constant,
the element $u$ constructed in Theorem \ref{HS}
is the unique horizontal section
with respect to $\nabla(\frac{d}{d\lam})$
near $\lam = 0$.
\end{Cor}

\pf
The Picard-Fuchs equation (\ref{PF})
satisfies the condition in Theorem \ref{HS}
with $A = \CC$, $B = \CC[[\lam]]$, $M = \MM_{dR}$,
$D = \lam \nabla(d/d\lam)$, $e = \lam$,
$g = F(\lam)$, and $\eta$ as defined in (\ref{eta}).
Hence $u$ is horizontal
with respect to $\lam \nabla(d/d\lam)$.
Thus it is horizontal 
with respect to $\nabla(\lam/d\lam)$.
The uniqueness follows
by the computation of the local monodromy near $\lam = 0$
(see \cite{HST}, Corollary 1.7).
\qed\\

\section{General properties
	and the method of Dwork and Katz}\label{Katz}

In this section,
we mainly consider the Dwork family
defined over the base $\BR = \FF_q$.
We define the Hasse invariant of the Dwork family
over characteristic $p$.
We derive the formula for the unit root
when the first slope of the Newton polygon of $V_t$
is zero via the crystalline approach.\\

\ni{\em (a) The Hasse invariant}\\

In the rest of this paper,
let $F(x)$ denote the hypergeometric {\em series}
(\ref{F_n_n-1})
\begin{eqnarray*}
F(x) &=& \shyperg{x} \\
	&=& \sum_{r=0}^{\infty} \frac{\pochf{1}{n+1}{r}
		\pochf{2}{n+1}{r} \cdots\pochf{n}{n+1}{r}}{(r!)^n}
		x^r \\
	&=& \sum_{r=0}^{\infty} \left( \prod_{i=1}^{n+1} \binom{ri}{i}\right)
		\left( \frac{x}{(n+1)^{n+1}} \right)^r.
\end{eqnarray*}
The last expression shows that
$F(x)$ is a formal power series
with coefficients in $\ZZ[\frac{1}{n+1}]$.
For any positive integer $s$,
let $F^{<s}(x)$ be the truncated polynomial of the series $F(x)$
up to degree $s-1$.
We will freely regard $F^{<s}(x)$
as a polynomial in any $\ZZ[\frac{1}{n+1}]$-algebra.
Similarly for the formal power series
$F(x)$.\\

\ni{\em Definition}.
Fix a prime $p$ with $(p, n+1) = 1$.
The function $H(x) = F^{<p}(x)$,
regarded as a polynomial over $\FF_p$,
is called the
{\em Hasse invariant}
of the Dwork family $V_t$
over characteristic $p$.\\

Let $A_m$ be the coefficient of $(X_1 \cdots X_{n+1})^m$
in $\Dw_t(X)^m$.
Thus $A_m$ is a polynomial in $t$ of degree $m$.\\

\begin{Lemma}\label{modp}
We have
\renewcommand{\labelenumi}{{\rm (\roman{enumi})}}
\begin{enumerate}
\item If $q = p^r$, then
$A_{q-1} \equiv (A_{p-1})^{(q-1)/(p-1)} \mod{p}$.
\item Let $\lam = t^{-(n+1)}$.
As polynomials in $t$, we have
$A_{p-1} \equiv t^{p-1} H(\lam) \mod{p}$.
\end{enumerate}
\end{Lemma}

\pf
(i)
We have
\begin{eqnarray*}
\Dw_t(X)^{q-1} &=& \left( \Dw_t(X)^{1+p+\cdots + p^{r-1}} \right)^{p-1} \\
	&\equiv& \left( \Dw_t(X) \Dw_{t^p}(X^p) \cdots
		\Dw_{t^{p^{r-1}}}(X^{p^{r-1}}) \right)^{p-1} \pmod{p}.
\end{eqnarray*}
Notice that
the coefficient $B_i$ of $(X_1 \cdots X_{n+1})^{p^i (p-1)}$
in $\Dw_{t^{p^i}}(X^{p^i})$
is congruent to $A_{p-1}^{p^i} \mod{p}$.
By an inspection of the terms in the product,
we have, in $\FF_p[t]$,
\begin{eqnarray*}
A_{q-1} &=& B_0 B_1 \cdots B_{r-1} \\
	&=& A_{p-1} A_{p-1}^p \cdots A_{p-1}^{p^{r-1}}.
\end{eqnarray*}

(ii)
We have
\begin{eqnarray*}
A_{p-1} &=& \sum_{i \geq 0} (-(n+1)t)^{p-1-(n+1)i}
		\binom{p-1}{i}\binom{p-1-i}{i} \cdots \binom{p-1-ni}{i} \\
	&=& (-(n+1)t)^{p-1} \sum_{i \geq 0}
		\left( (-(n+1))^{-(n+1)i} \prod_{k=0}^n \binom{p-1-ki}{i} \right)
			\lam^i \\
	&\equiv& t^{p-1} \sum_{i=0}^{\left[\frac{p-1}{n+1}\right]}
		\left( \prod_{k=0}^{(n+1)i} \frac{k}{n+1} \right)
			\frac{\lam^i}{(i!)^{n+1}} \pmod{p} \\
	&\equiv& t^{p-1} \sum_{i=0}^{p-1}
		\left( \prod_{k=0}^{(n+1)i} \frac{k}{n+1} \right)
			\frac{\lam^i}{(i!)^{n+1}} \pmod{p} \\
	&=& t^{p-1} H(\lam).
\end{eqnarray*}
\qed

\begin{Thm}\label{Hasse}
Let $\FF_q$ be a finite field of $q$ elements
with $q$ a power of $p$
and $(p, n+1) = 1$.
Let $H(x) = F^{<p}(x)$
be the Hasse invariant of the Dwork family.
Let $t \in \FF_q, t \neq 0, t^{n+1} \neq 1$.
Let $\lam = t^{-(n+1)}$.
Then the first slope of the Newton polygon of $V_t$ is zero
if and only if
$H(\lam) \neq  0$.
In this case,
if $\pi_t$ is the unit root of $V_t$,
then
\[ \pi_t \equiv H(\lam)^{(q-1)/(p-1)} \pmod{p}. \]
\end{Thm}

\pf
(a)
The method here is similar to the case of counting points
of Legendre family of elliptic curves.
Let
\[ N_t = \# \{ x \in \PP^n(\FF_q) | x \in V_t (\FF_q) \}, \quad
\text{and}\quad N'_t = \# \{ x \in \FF_q^{n+1} | \Dw_t(x) = 0 \}. \]
Then
\[ N_t = \frac{N'_t - 1}{q-1} \equiv 1 - N'_t \pmod{q}. \]
Consider the zeta-function of $V_t$
\begin{eqnarray*}
Z(V_t, T) &:=& \exp \left(N_t T + \mathcal{O}\left(T^2 \right) \right) \\
	&\equiv& 1 + N_t T \pmod{T^2}.
\end{eqnarray*}
Let
\[ \det\left( 1 - T\ {\rm Frob}^* \big|
	H^{n-1}_{\text{\'et}}(V_t \tensor \bar\FF_q, \QQ_\ell) \right)
		=1-aT + \mathcal{O}(T^2) \]
be the reciprocal characteristic polynomial
of the geometric Frobenius ${\rm Frob}^*$
on the middle cohomology of $V_t$ ($\ell \neq p$).
Then by the Weil conjecture
\begin{eqnarray*}
Z(V_t, T) &=& \frac{\left(1-aT+ \mathcal{O}(T^2)\right)^{(-1)^n}}
		{(1-T)(1-qT^2) \cdots (1-q^nT^{2n})} \\
	&\equiv& \left( 1- (-1)^n aT \right)(1+T) \pmod{T^2} \\
	&\equiv& 1 + \left( 1- (-1)^n a \right) T \pmod{T^2}.
\end{eqnarray*}
Thus $a \equiv (-1)^n N_t' \mod{q}$,
and the theorem is equivalent to say that
for $\lam \neq 0, 1$,
the congruence $N_t' \equiv 0 \mod{p}$ holds
if and only if
$H(\lam) = 0$ in $\FF_q$.

By Warning's method,
\[ N'_t = \sum_{x \in \FF_q^{n+1}} \left( 1 - \Dw_t(x)^{q-1} \right). \]
Notice that
\[ \sum_{z \in \FF_q} z^r \equiv \left\{ \begin{array}{cl}
	-1 & \text{if $(q-1) | r$} \\
	0 & \text{otherwise}. \end{array} \right. \]
Thus by an inspection of the terms
in the expansion of $\Dw^{q-1}$,
\[ (-1)^n \sum_{x \in \FF_q^{n+1}} \left( 1 - \Dw_t(x)^{q-1} \right)
	\equiv \text{the coeff. $A_{q-1}$ of $(X_1 \cdots X_{n+1})^{q-1}$
		in $\Dw^{q-1}$} \pmod{p}. \]
Therefore by Lemma \ref{modp},
the congruence $N_t' \equiv 0 \mod{p}$ holds
if and only if
\[ 0 = \left( t^{p-1} H(\lam) \right)^{(q-1)/(p-1)} = H(\lam)^{(q-1)/(p-1)}, \]
since $t^{q-1} = 1$.

(b)
Suppose now $H(\lam) \neq 0$,
then we have
\[ \pi_t \equiv a \equiv A_{q-1} \equiv H(\lam)^{(q-1)/(p-1)}
	\pmod{p} \]
by the above calculation.
\qed\\

\ni{\em Remark.}
The proof of the theorem also shows that
over $\FF_q$,
the first slope of the variety $V_0$ in the family at $t=0$
is zero
if and only if
$H(x)$ is strictly of degree $(p-1)/(n+1)$
as an element in $\FF_q[x]$.
This condition is equivalent
to the congruence
$p \equiv 1 \mod{(n+1)}$.
The Newton polygon of $V_0$
can be determined
by looking at the splitting type of the prime $p$
in the cyclotomic field $\QQ(\mu_{n+1})$
(see \cite{GY}, Proposition 3.8).\\

\ni{\em (b) The formula for the unit root}\\

\begin{Thm}\label{unitroot}
Let $p$ be a prime with $(p,n+1)=1$.
Let
\[ f(x) = \frac{F(x)}{F(x^p)} \]
as a formal power series in $\ZZ_p[[x]]$.
Then the following assertions
hold true.
\renewcommand{\labelenumi}{{\rm (\arabic{enumi})}}
\begin{enumerate}
\item
$f(x)$ is in fact an element
in the $p$-adic completion of
$\ZZ_p[x, (x(1-x)H(x))^{-1}]$.
\item
Let $q = p^r$
and $t \in \FF_q$, $t \neq 0, t^{n+1} \neq 1$.
Put $\lam = t^{-(n+1)}$.
Let $\hat{\lam}$ be the Teichm\"uller lifting of $\lam$.
If $H(\lam) \neq 0$,
then
\[ \pi_{\lam} =
	f(\hat{\lam}) f(\hat{\lam}^p) \cdots f(\hat{\lam}^{p^{r-1}}) \]
is the unique unit root of $V_t$.
\end{enumerate}
\end{Thm}

\pf
(1)
Let
\[ R = W(\bar{\FF}_p)[t, (t(1-t^{n+1})H(\lam))^{-1}]. \]
Let $S = \Spec R$,
$S_0 = \Spec R/p R$,
and $S^{\infty} = \Spec \projlim R/p^n R$,
where the projective limit runs over all $n>0$.
Choose a lifted Frobenius $\sigma$ on $S^{\infty}$
by taking $\sigma(t) = t^p$.
The pointwise defined $\MM_{cris, t}$ in \S $2(b)$
forms a Hodge $F$-crystal $\MM_{cris}$ of rank $n$
on $S^{\infty}$.
As a sheaf of modules,
$\MM_{cris}$ is isomorphic
to the subsheaf $\MM_{dR}$
of the relative de Rham cohomology
$R^{n-1}v^{\infty}_* \Omega_{\mathcal{V}^{\infty}/S^{\infty}}^{\bullet}$
of the Dwork family
$\mathcal{V}^{\infty}$ over $S^{\infty}$
(see \cite{Katz_Dwork}, \S 7).
Let $\phi$ be the absolute Frobenius on $\MM_{cris}$
with respect to $\sigma$.

Over each geometric point $t$ of $S_0$,
the Newton polygon of the crystal $\MM_{cris, t}$
begins with a segment of slope zero of length 1
(Theorem \ref{Hasse}).
On the other hand, the absolute Frobenius $\phi$
on the Hodge filtration
${\rm Fil}^1 \rc \MM_{cris}$
satisfies
\[ \phi \left( {\rm Fil}^1 \right) \rc p \cdot \MM_{cris} \]
since $\MM_{cris}$ is from geometry
(see \cite{Katz_Dwork}, \S 7).
Therefore
the unit root sub-crystal $U$ of $\MM_{cris}$
is generated over $W(\bar{\FF}_q)$
by horizontal sections (op.cit., 4.1.2).
Notice that
the Picard-Fuchs equation (\ref{PF})
has a unique power series solution $F(\lam)$
in $W(\bar{\FF}_p)[[\lam]]$
with constant term 1.
Therefore by Corollary \ref{HS_Dwork},
the crystal $U$ is generated by $u$ defined in Theorem \ref{HS}.
Thus
the series $f(\lam)$ is an element in $\projlim R/p^n R$
(op.cit., 4.1.9).
Since $f(\lam)$ depends only on $\lam$
and its coefficients are $p$-adic integers,
the assertion follows.

(2)
Since $U$ is generated over $W(\bar{\FF}_q)$ by $u$,
there exists an $c \in W(\bar{\FF}_p)$
such that $cu$ is fixed by the Frobenius
and $\ve f(\hat{\lam})$
(with $\ve = c^{1-\sigma}$)
represents the unit root
of the absolute Frobenius on $\MM_{cris,\lam}$
with respect to some bases
over $W(\FF_q)$
(cf. \cite{Katz_Dwork}, \S 8).
Thus we obtain
\[ \pi_{\lam} = \ve^{1 + \sigma + \cdots + \sigma^{p^{r-1}}}
	f(\hat{\lam}) f(\hat{\lam}^p) \cdots f(\hat{\lam}^{p^{r-1}}). \]
Let
\[ \ve' := \ve^{1 + \sigma + \cdots + \sigma^{p^{r-1}}} \]
be the constant term above.
We now ought to show that
$\ve' = 1$.

We apply Lemma (6.2) in \cite{Dwork}
(cf. op.cit., \S $6(j)$ for $n = 3$).
In our case,
the nilpotent part
of the local monodromy (over characteristic 0)
near $t = \infty$
with respect to some bases
is given by the matrix
\[ \mathcal{N} = \left( \begin{array}{ccccc}
	0 & 1 & \ldots & 0 & 0 \\
	0 & 0 & \ldots & 0 & 0 \\
	\vdots & \vdots &\ddots & \vdots & \vdots \\
	0 & 0 & \ldots & 0 & 1 \\
	0 & 0 & \ldots & 0 & 0 \end{array} \right), \]
which is nilpotent of order $n$ (\cite{HST}, Corollary 1.7).
By loc.cit.,
there exist scalar $(n \times n)$ matrices
$A, \Theta$
such that
\begin{equation}\label{c1}
A \mathcal{N} = q \mathcal{N} A
\end{equation}
\begin{equation}\label{c2}
\Theta N^t + N \Theta = 0
\end{equation}
\begin{equation}\label{c3}
q^{n-1} \Theta = A \Theta A^t
\end{equation}
and such
that $\ve'$ is the $(n,n)$-th component of $A$.
Under the conditions
(\ref{c1}), (\ref{c2}) and (\ref{c3}),
one finds that
\[ A = \left( \begin{array}{cccc}
	q^{n-1} \ve' & * & * & * \\
	0 & q^{n-2} \ve' & * & * \\
	\vdots & \vdots & \ddots & * \\
	0 & 0 & \ldots & \ve' \end{array} \right) \]
and $\ve' = \pm 1$.
By op.cit. Lemma (3.4)(i),
\[ f(x) \equiv H(x) \pmod{p}. \] 
By the second part of Theorem \ref{Hasse},
we must then have $\ve' = 1$.
\qed\\

\begin{Prop}
Let
\[ F_i(x) = \frac{d^i}{dx^i} F(x) \]
be the $i$-th derivative of $F(x)$.
Then
the series
\[ f_i(x) = \frac{F_i(x)}{F(x)} \]
are in fact elements
in the $p$-adic completion of
$\ZZ_p[x, (x(1-x)H(x))^{-1}]$
for all $i \geq 0$.
\end{Prop}

\pf
We keep the notations as in the proof of Theorem \ref{unitroot}.
Write
\[ u = \sum_{i=0}^{n-1} C_i \eta^{(n-1-i)}. \]
Since the unit crystal $U$ is generated by $u$,
it follows that
$C_i/C_0$ are elements in $\projlim R/p^n R$
for all $0 \leq i \leq (n-1)$
(\cite{Katz_Dwork}, 4.1.9).
By the explicit description of $u$
(see Remark after Theorem \ref{HS}),
we see inductively that
$f_i(\lam)$ is in $\projlim R/p^n R$.
Since $f_i(\lam)$ depends only on $\lam$
and has $p$-adic integral coefficients,
the assertion follows for $0 \leq i \leq n-1$.

Since $F(x)$ is a solution to (\ref{PF}),
the higher derivatives $F_i(x)$
can be written as a
$\ZZ_p[x, (x(1-x))^{-1}]$-combination
of $\{ F_i(x) \}_{0 \leq i \leq n-1}$.
Thus the assertion also holds for all $i \geq 0$.
\qed\\

\section{Method of Stienstra and Beukers}\label{SB}

In this section,
we study the unit root of $V_t$
form the point of view of formal groups.
Recall that
$A_m$ is the coefficient of $(X_1 \cdots X_{n+1})^m$
in $\Dw_t(X)^m$.\\

\ni{\em (a) The formal group laws}\\

\begin{Lemma}\label{coeff}
Let $\lam = t^{-(n+1)}$.
As polynomials in $t$,
we have
\[ A_m = (-(n+1)t)^m
\hyperg{n+1}{n}{\frac{-m}{n+1}}{\frac{-m+1}{n+1}}{\frac{-m+n}{n+1}}{1}{1}{1}
	{\lam}. \]
\end{Lemma}

\pf
We have
\begin{eqnarray*}
A_m &=& \sum_{r \geq 0} \left( \prod_{i=0}^n \binom{m-ri}{r} \right)
	(-(n+1)t)^{m-r(n+1)} \\
&=& (-(n+1)t)^m \sum_{r \geq 0} \binom{m}{(n+1)r}
	\frac{((n+1)r)!}{(-(n+1))^{(n+1)r}(r!)^{n+1}}
		\left( \frac{1}{t^{n+1}} \right)^r \\
&=& (-(n+1)t)^m
\hyperg{n+1}{n}{\frac{-m}{n+1}}{\frac{-m+1}{n+1}}{\frac{-m+n}{n+1}}{1}{1}{1}
	{\lam}.
\end{eqnarray*}
\qed

\begin{Prop}\label{law}
Consider the family defined by equation (\ref{Dwork})
over a noetherian ring $\BR$, which is flat over $\ZZ$.
Let $t \in \BR$ and $\lam = t^{-(n+1)}$.
The formal group $H^{n-1}_{et}(V_t, \hat{\GG}_m)$
associated to $V_t$
can be realized as the formal group law $G_t$ over $\BR$
with logarithm $l(\tau) \in \QQ \tensor_{\ZZ} \BR[[\tau]]$ given by
\[ l(\tau) = \sum_{m=0}^{\infty} (-(n+1)t)^m
\hyperg{n+1}{n}{\frac{-m}{n+1}}{\frac{-m+1}{n+1}}{\frac{-m+n}{n+1}}{1}{1}{1}
	{\lam}
	\frac{\tau^{m+1}}{m+1}. \]
\end{Prop}

\pf
The formal group law with logarithm
\[ \sum_{m=0}^{\infty} A_m \frac{\tau^{m+1}}{m+1} \]
realizes the formal group associated to $V_t$
(\cite{S}, Theorem 1).
Thus the statement follows by Lemma \ref{coeff}.
\qed\\

\ni{\em Remark.}
By varying $t$,
the formal groups $G_t$ above
can be put together
as a flat family of formal groups
over the base space $t \in \PP^1$.
When $t = \infty$,
the formal group $G_{\infty}$
associated to
\[ V_{\infty}: -(n+1)X_1 X_2 \cdots X_{n+1} \]
is given by the logarithm
\[ l(\tau) = \sum_{m=0}^{\infty} (-(n+1))^m \frac{\tau^{m+1}}{m+1}. \]
Thus
\begin{eqnarray*}
G_{\infty}(x,y) &=& l^{-1} \left( l(x) + l(y) \right) \\
	&=& x + y + (n+1)xy.
\end{eqnarray*}
Via the transformation $x \mapsto (n+1)x$,
the group $G_{\infty}$
is isomorphic over $\BR$
to the standard multiplicative formal group
\[ \hat{\GG}_m(x,y) = x + y + xy. \]
If the residue field of $\BR$ at a closed point
is a finite field $\FF_q$,
then the Frobenius endomorphism
on the reduction of $G_{\infty}$
at that point
acts as multiplication by $q$.\\

\ni{\em (b) Proofs of main results in $\S \ref{Katz}$ via formal groups}\\

In what follow,
let $\RR$ be the $p$-adic completion
of the ring $\ZZ_p [x, (x(1-x)H(x))^{-1}]$.
Let $\sigma$ be the endomorphism of $\RR$
extending the Frobenius on the constants
and with $\sigma(x) = x^p$.
For $a \in \RR$,
we write $a^{\sigma} = \sigma(a)$.\\

To facilitate the discussion,
we write
\[ F_{m,s}(x) = F^{< mp^s} (x) \]
for the truncated hypergeometric series
up to degree $mp^s - 1$.
Recall that $H(x) = F_{1,1}(x)$
is the Hasse invariant.
Let
\[ G_{\mu,s}(x) =
	\hyperg{n+1}{n}{\frac{-\mu p^s+1}{n+1}}{\frac{-\mu p^s+2}{n+1}}
	{\frac{-\mu p^s+n+1}{n+1}}
		{1}{1}{1}{x}. \]
Note that it is a polynomial of degree
$\left[ \frac{\mu p^s-1}{n+1} \right]$.
Here $[z]$ denotes
the least integer function.
Also write (with $\lam = t^{-(n+1)}$)
\[ G'_{\mu,s}(t) = (-(n+1)t)^{\mu p^s -1} \cdot G_{\mu,s}(\lam) \]
for the coefficient of $\frac{\tau^{\mu p^s}}{\mu p^s}$
in the logarithm $l(\tau)$ of $G_t$ in Proposition \ref{law}.
We regard $G'_{\mu,s}(t)$ as a polynomial in $t$.\\

Let $\CC_p$ be the $p$-adic completion
of an algebraic closure of $\QQ_p$.
Take an non-archimedean norm
$| \cdot |$ on $\CC_p$.\\

\begin{Lemma}\label{cong}
Regarding $H$ and $G_{\mu,s}$ as elements in $\mathcal{R}$,
we have
\renewcommand{\labelenumi}{{\rm (\roman{enumi})}}
\begin{enumerate}
\item
If $1 \leq \mu \leq n+1$,
then $G_{\mu,1} \equiv H \mod{p}$.
\item
$G_{\mu,s+1} \equiv G_{\mu,0}^{\sigma^{s+1}} \cdot
	H^{1+\sigma + \cdots \sigma^s} \mod{p}$.
\item
There exists an element $g \in \RR$ such that
$G_{\mu,s+1} \equiv g \cdot G_{\mu,s}^{\sigma} \mod{p^{s+1}}$
for all $\mu, s \geq 0$.
\item
For any $\lam \in \CC_p$,
if $|H(\lam)| = 1$,
then $|G_{\mu,s}(\lam)| = 1$.
\end{enumerate}
\end{Lemma}

\pf
(i) is obvious by the observation
that $G_{\mu,1}(\lam)$ has degree $[\frac{\mu p-1}{n+1}] < p$
if $\mu \leq n+1$.

(ii) and (iii) are direct consequences of results in \cite{SB}.
Let $\mathcal{S}$ be the $p$-adic completion
of the ring $\ZZ_p[t, (t(1-t^{n+1})H(\lam))^{-1}]$.
Let $\sigma$ be the Frobenius on $\mathcal{S}$
with $t^{\sigma} = t^p$.
Notice that by definition,
$t^{p-1}H(\lam)$ is invertible in $\mathcal{S}$.
Thus the reduction of the formal group law
defined in Proposition \ref{law}
to any point of $\mathcal{S}$ of characteristic $p$
is of multiplicative type
(see \cite{SB}, Theorem (A.8)(v)).
This implies (loc.cit.) that
there exists an element $g' \in \mathcal{S}$
such that
\[ G'_{\mu,s+1} \equiv g' \cdot (G'_{\mu,s})^{\sigma} \pmod{p^{s+1}}. \]
Thus
\[ (-(n+1))^{\mu p^s(p-1)} t^{p-1} G_{\mu,s+1}
	\equiv g' \cdot G_{\mu,s}^{\sigma} \pmod{p^{s+1}}. \]
Since $(-(n+1))^{\mu p^s(p-1)} \equiv 1 \mod{p^{s+1}}$,
we have
\[ G_{\mu,s+1} \equiv \frac{g'}{t^{p-1}} G_{\mu,s}^{\sigma}
	\pmod{p^{s+1}}. \]
Let $g = g' / t^{p-1}$.
Then $g$ depends only on $\lam$
and hence it is obvious that indeed $g(x) \in \RR$.

Since $g \equiv G_{1,1} \equiv H \mod{p}$,
inductively we get (ii).

For (iv), assume that $|H(\lam)| = 1$.
If $\mu < n+1$,
then $G_{\mu,0} = 1$.
Thus $|G_{\mu,s}(\lam)| = 1$ by (ii).
In general,
we can choose some $\ve$ such that
$G_{\mu,s+\ve} = G_{\mu',s'}$ with $0 \leq \mu' < n+1$.
Since
\[ G_{\mu,s+\ve} \equiv G_{\mu,s}^{\sigma^{\ve}} \cdot
	H^{1 + \sigma + \cdots + \sigma^{\ve-1}} \pmod{p}, \]
the equality $|G_{\mu',s'}(\lam)| = 1$ implies $|G_{\mu,s}(\lam)| = 1$.
\qed\\

\ni{\em The second proof of Theorem \ref{Hasse}.}\\

Take a lifting $\hat{t} \in W(\FF_q)$ of $t$
and let $\hat{\lam} = \hat{t}^{-(n+1)}$.
Then $G_{\hat{t}}$
constructed in Proposition \ref{law}
is a formal group over $W(\FF_q)$
whose reduction to $\FF_q$ is the formal group $G_t$
associated to $V_t$.
The group $G_t$
is of height one
if and only if
the coefficient of $\tau^{p}/p$
in the logarithm $l(\tau)$ of $G_{\hat{t}}$
is invertible in $W(\FF_q)$.
By Lemma \ref{cong}(i),
\[ \hyperg{n+1}{n}{\frac{-m}{n+1}}{\frac{-m+1}{n+1}}{\frac{-m+n}{n+1}}
		{1}{1}{1}{\hat{\lam}} = G_{1,1}(\hat{\lam})
	\equiv H(\lam) \pmod{p} \]
and hence the first assertion follows.

The second assertion follows from Lemma \ref{cong}(i).
Notice that the remark after Theorem \ref{Hasse}
also follows easily by the same argument.
\qed\\

\ni{\em The second proof of Theorem \ref{unitroot}}.\\

(1)
Observe that as $s \to \infty$,
the elements $G_{\mu,s}(x)$ converge to $F(x)$ $p$-adically.
Therefore we see that
the element $g \in \RR$ in Lemma \ref{cong}(iii)
must be
\begin{equation}\label{g=f}
\frac{F(x)}{F(x)^{\sigma}} = \frac{F(x)}{F(x^p)} = f(x).
\end{equation}

(2)
This is an exercise on formal group theory
and all we need is already in \cite{SB}.
Let $\hat{t}$ be the Teichm\"uller lifting of $t$.
Let $a = \hat{t}^{p-1} f(\hat{\lam})$
and
\[ a_s = a^{1 + \sigma + \cdots \sigma^{s-1}}. \]
Consider the formal group law $G'$ over $W(\FF_q)$
with logarithm
\[ l'(\tau) = \tau + \sum_{s \geq 1} a_s \frac{\tau^{p^s}}{p^s}. \]
Notice that
the coefficients of $l'(\tau)$
satisfy $a_{s+1}/a_s^{\sigma} = a$
for any $s \geq 0$.
Thus the formal groups $G$ and $G'$
are strictly isomorphic to each other
over $W(\FF_q)$
(\cite{SB}, Theorems (A.8) and (A.9)).

On the other hand,
the group $G'$ is isomorphic to a formal group
whose Cartier module
(which is of rank one over $W(\FF_q)$)
has a basis $\omega$
with the Frobenius acting as $\omega \mapsto a \omega$
(op.cit. (A.13)).
Thus the $p$-adic unit
\[ f(\hat{\lam}) f(\hat{\lam}^p) \cdots f(\hat{\lam}^{p^{r-1}})
	= a a^{\sigma} \cdots a^{\sigma^{r-1}} \]
equals to the Frobenius endomorphism
of the Cartier module of $G$.
Hence it is the unique eigenvalue
of the geometric Frobenius endomorphism
on the middle cohomology of $V_t$.
\qed

\ni{\em Remark.}
The formula for the unit root
also make sense when $\lam = 0$
if one consider the variation
of the Frobenius endomorphism
on the Cartier module
of the formal groups $G_t$ associated to $V_t$.
See the remark after Proposition \ref{law}.\\

\ni{\em (c) Dwork's congruences}\\

Here we remark some congruent relations.\\

Write $F(x) = \sum B(i) x^i$.
Combine Lemma \ref{cong} (i) and (ii),
we find
\[ F(x) \equiv F^{< p}(x) F(x^p) \pmod{p}. \]
Comparing the coefficients on both sides,
this implies for $0 \leq c < p$,
we have
\[ B(c+p) \equiv B(c) B(1) \pmod{p}, \]
which is a tiny special case of
\cite{Dwork}, \S 1, Corollary 2.\\

In \cite{Dwork}, Lemma (3.4),
one finds the congruences
\[ F_{m,s+1}(x) \cdot F(x^p)
	\equiv F_{m,s}(x^p) \cdot F(x) \pmod{p^{s+1}} \]
for any integers $m, s \geq 0$.
Thus by Lemma \ref{cong} and (\ref{g=f}),
we have that
\begin{equation}\label{twohyperg}
F_{m,s+1}(x) \cdot G_{\mu,s}(x^p)
	\equiv F_{m,s}(x^p) \cdot G_{\mu,s+1}(x) \pmod{p^{s+1}}
\end{equation}
for any non-negative integers $m, \mu, s$
as polynomials in $\ZZ[\frac{1}{n+1}][x]$.
Consider the formal group laws
$J_{\lam}$ over $\ZZ_p[\lam]$
and $J_{t}$ over $\ZZ_p[t, t^{-1}]$
with logarithms
\[ j_{\lam} = \sum F_{m,s}(\lam) \frac{\tau^{mp^s}}{mp^s}
\quad \text{and}\quad
j_{t} = \sum t^{mp^s -1} F_{m,s}(\lam) \frac{\tau^{mp^s}}{mp^s}, \]
respectively.
Then the congruences (\ref{twohyperg})
imply that
$J_t$ is strictly isomorphic to $G_t$
over $\ZZ_p[t, t^{-1}]$.
It would be interesting to see
if one can find some geometry
behind these formal groups
and its relation to the Dwork family $V_t$.

\end{document}